\documentstyle[11pt]{article}

\begin{document}

\input amssym.def
\input  amssym.tex

\title{Natural Internal Forcing Schemata extending ZFC. A Crack in the
Armor surrounding $CH?$}

\author{Garvin Melles\thanks{Would like to thank Ehud Hrushovski
for supporting him with funds from NSF Grant DMS 8959511} \\Hebrew University of Jerusalem}

\newtheorem{theorem}{Theorem}[section]
\newtheorem{defi}[theorem]{Definition}
\newtheorem{lemma}[theorem]{Lemma}
\newtheorem{coro}[theorem]{Corollary}
\newtheorem{conj}{Conjecture}

\newcommand{\proof}{{\sc proof} \hspace{0.1in}}
\newcommand{\iopp}{\stackrel{|}{\smile}}
\newcommand{\niopp}{\not\!\!{\stackrel{|}{\smile}}}
\newcommand{\noindep}[1]{\mathop{\niopp}\limits_{\textstyle{#1}}}
\newcommand{\indep}[1]{\mathop{\iopp}\limits_{\textstyle{#1}}}
\newcommand{\sub}{\subseteq}

\mathsurround=.1cm
\maketitle

Mathematicians are one over on the physicists in that they already have a
unified theory of mathematics, namely set theory. Unfortunately the
plethora of independence results since the invention of forcing has
taken away some of the luster of set theory in the eyes of many
mathematicians. Will man's knowledge of mathematical truth
be forever limited to those theorems derivable from the standard
axioms of set theory, $ZFC?$ This author does not think so, and in fact
he feels there is a schema concerning non-constructible sets
which is a very natural candidate for being  considered as part of the
axioms of set theory. To
understand the motivation why, let us take a very short look back at
the history of the development of mathematics. Mathematics began
with the study of mathematical objects very physical and concrete in
nature and has progressed to the study of
things completely imaginary and abstract. Most mathematicians now accept these objects
as as mathematically legitimate as any of their more concrete
counterparts. It is enough that these objects are consistently
imaginable, i.e., exist in the world of
set theory.  Applying the same intuition to set theory itself, we
should accept as sets as many that we can whose existence are
consistent with $ZFC.$ Of course this is only a vague notion, but
knowledge of set theory so far, namely of the existence of $L$
provides a good starting point. What sets can we consistently imagine
beyond $L?$ 
Since by forcing one can prove the
consistency of $ZFC$ with the existence of non-constructible sets and as $L$ is
absolute, with these forcing extensions of $L$ you have consistently imagined more sets in
a way which satisfies the vague notion mentioned above. The problem
is which forcing extensions should you consider as part of the
universe? But there is no problem, because if you prove the
consistency of the existence of some $L$ generic subset of a partially
ordered set $P\in L$ with $ZFC,$ then $P$ must be describable
and we can easily prove the consistency of $ZFC$ with the existence of
$L$ generic subsets of $P$ for every $P$ definible in $L.$ Namely, the
axiom schema $IFS_L$ (For internal forcing schema over $L\!$) defined
below is consistent with $ZFC.$

\begin{defi}
$IFS_L$ is the axiom schema which says for every formula $\phi(x),$ if
$L\models$ there is a unique partial order $P$ such that $\phi(P),$
then there is a $L$ generic subset of $P$ in the universe $V.$
\end{defi}

\noindent $IFS_L$ is a natural closure condition on a universe
of set theory. Given a class model of $ZFC$ which has no inner class model of the
form $L[G]$ for some partial order $P$ definable in $L,$ we can (by
forcing) consistently imagine expanding the model to include such a
class. Conversely, no class  model of $ZFC+IFS_L$ can
be contained in a class model of $ZFC$ which does not satisfy $IFS_L.$

\begin{theorem}
If there is a sequence $\langle M_n\mid n<\omega\rangle$ of transitive models with
$M_n\models ZFC_n$ where $ZFC=\bigcup\limits_{n\in\omega}ZFC_n,$ then $Con(ZFC+IFS_L)$
\end{theorem}
\proof By the compactness theorem and forcing. 

\begin{theorem}\label{IFSLCBLG}
If $V$ is a model of $ZFC,$ then $V\models IFS_L$ if and only if
$V\models$ every set definable in $L$ is countable. 
\end{theorem}
\proof Certainly if every set definable in $L$ is countable, then
every partially ordered set definable in $L$ is countable, 
so therefore is the set of
dense subsets of $P$ in $L$ countable and so
$P$ has generic subsets over $L$ in the universe. In the other
direction, if $s$ is a set definable in $L,$ then so is the partially
ordered set consisting of maps from distinct finite subsets of $s$ to distinct finite
subsets of $\omega,$ so a $L$ generic subset over the partial ordering
is a witness to $\big|s\big|=\omega.$

\vspace{.1in}

\noindent Perhaps $IFS_L$ is not surprising since $ZFC+0^{\#}\ \vdash\ IFS_L.$
But the same reasoning as led to $IFS_L$ leads to the following
stronger schema, $IFS_{Ab\,L[r]}$ (For internal forcing schema for absolute
class models of $ZFC$ constructible over an absolutely definable real) which implies that if
$0^{\#}$ exists, then all sets definable in $L[0^{\#}]$ are countable.

\begin{defi}
A subset $r$ of $\omega$ is said to be absolutely definable if for
some $\Pi_1$ formula $\theta(x),$

\begin{enumerate}
\item $V\models \theta(r)$
\item $ZFC\vdash \exists\,x\theta(x)\ \rightarrow\
\exists!\,x\theta(x)$

\end{enumerate}  

\end{defi}

\begin{defi}
$IFS_{Ab\,L[r]}$ is the axiom schema of set theory which says if $r$
is an absolutely definable real then all definable elements of $L[r]$
are countable (equivalently, every partial order $P$ definable in
$L[r]$ has an $L[r]$ generic subset.)
\end{defi}

\noindent The following theorem is a formal justification of $IFS_{Ab\,L[r]}.$

\begin{theorem}\label{CONIFSr}
Suppose $V$ is a countable transitive model of $ZFC$ and let $\big\{\theta_i(x)\mid
i<\omega\big\}$ be the list of all formulas defining absolute reals
such that $V\models
\bigwedge\limits_{i<\omega}\exists\,x\theta_i(x).$  Suppose that the
supremum of the ordinals definable in $V$ is in $V.$ Then there is a
countable transitive extension $V'$ of $V$ with the same ordinals such that
$$V'\models ZFC+IFS_{Ab\,L[r]}+\bigwedge\limits_{i<\omega}\exists\,x\theta_i(x)$$
\end{theorem}
\proof Let $\alpha^*$ be the sup of all the ordinals definable in $L.$
Let $P$ be the set of finite partial one to one functions from
$\alpha^*$ to $\omega.$ Let $V'=V[G]$ where $G$ is a $V$ generic
subset of $P.$ 
To finish the proof it is enough to prove the following claim.

\vspace{.1in}

\noindent Claim: If $\psi(x)$ defines a real in $M[G]$ then it
is in $M.$\\
\proof Since $P$ is separative, if $p\in P$ and $\pi$ is an
automorphism of $P,$ then for every formula $\varphi(v_1,\ldots,v_n)$
and names $x_1,\ldots,x_n$ 
$$*\ \ \ \ p\Vdash \varphi(x_1,\ldots,x_n)\ \hbox{ iff }\ \pi p\Vdash
\varphi(\pi x_1,\ldots,\pi x_n)$$
Let $\varphi(x)=\exists\,Y(\psi(Y)\ \wedge\ x\in Y).$ Let $n\in
\omega.$ If for no $p\in P$ does $p\Vdash
\big|\big|\varphi(\check n)\big|\big|$ then
$\big|\big|\varphi(\check n)\big|\big|=0.$ So let $p\in P$ such that $p\Vdash
\big|\big|\varphi(\check n)\big|\big|.$ By $*$ if $\pi$ is an automorphism
of $P$ then  $\pi p\Vdash
\big|\big|\varphi(\check n)\big|\big|.$ Let $\pi$ be a permutation of
$\omega.$ $\pi$ induces a permutation of $P$ by letting for $p\in P,$
$dom\,\pi p=dom\,p$ and letting $\pi p(\alpha)=\pi(p(\alpha)).$ By
letting $\pi$ vary over the permutations of $\omega$ it follows that
$\big|\big|\varphi(\check n)\big|\big|=1.$ Let $\dot r$ be the name
with domain $\big\{\check n\mid n<\omega\big\}$ and such that
$$\dot r(\check n)=\big|\big|\varphi(\check n)\big|\big|$$ 
$i_G(\dot r)=r,$ but then
$r=\big\{n\mid \big|\big|\varphi(\check n)\big|\big|=1\big\}$ which
means it is in $M.$

\begin{coro}
$ZFC+IFS_{Ab\,L[r]}+$'there are no absolutely definable
non-constructible reals' is
consistent. (Relative to the assumption of a countable transitive
model of $L$ with its definable ordinals having a supremum in the
model.)
\end{coro}

\vspace{.1in}

\noindent Since classes of the form $L[r]$
are absolute if $r$ is an absolutely definable real, they provide
reference points from which to measure the size of the universe.
We can extend the schema $IFS_{Ab\,L[r]}$ by exploiting the similarity
between a
class such as $L({\Bbb R})$ and a class of the form $L[r]$ where $r$
is an absolutely definable real. We can
argue that if $P$ is a partial order definable in $L({\Bbb R}),$  and
if a $V$ generic subset of $P$ cannot add any reals to $V,$ then an
$L({\Bbb R})$ generic subset of $P$ should exist in $V.$ $L({\Bbb R})$ is
concrete in the sense the interpretation of $L({\Bbb R})$ is absolute
in any class model containing  ${\Bbb R}$, and thereby
like classes of the form $L[r]$ where $r$ is an absolutely definable
real, $L({\Bbb R})$ provides a reference point from which to measure
the size of the universe.   
This leads to the following  natural strengthening of $IFS_L$ and 
$IFS_{Ab\,L[r]}.$ 



\begin{defi}
$x\in V$ is said to be weakly absolutely definable of the form
$V_{\alpha}$ if for some
formula $\psi(v)$ which provably defines an ordinal and which is
provably $\Delta_1$ from $ZF,$ 
$$V\models \exists!\alpha\Big(\psi(\alpha)\ \wedge\ \forall y(y\in x\
\leftrightarrow\ \rho(x)\leq\alpha)\Big)$$ 
Let $\theta(x)$ denote $\exists!\alpha\Big(\psi(\alpha)\ \wedge\ \forall y(y\in x\
\leftrightarrow\ \rho(x)\leq\alpha)\Big)$ and let $ZF_{\theta}$ be a
finite part of $ZF$ which proves $\psi(v)$ is $\Delta_1$ and proves
$\psi(v)$ defines an ordinal. $\theta(x)$
is said to define a weakly absolutely definable set of the form
$V_{\alpha}.$ ($\rho(x)$ denotes the foundation rank.)
\end{defi}


\begin{defi}
$IFS_{W\!Ab\,L(V_{\alpha})}$ is the axiom schema of set theory which says for every
weakly absolutely definable set of the form $V_{\alpha}$ for every
partial order 
$P$ definable in
$L(V_{\alpha}),$ if
$$\big|\big|V_{\alpha}^{V[G]}=V_{\alpha}^{V}\big|\big|^{(r.o.P)^V}=1$$
then there exists an $L(V_{\alpha})$ generic subset $G$ of $P.$
\end{defi}

\begin{theorem}
If there is a sequence $\langle M_n\mid n<\omega\rangle$ of transitive models with
$M_n\models ZFC_n$ where $ZFC=\bigcup\limits_{n\in\omega}ZFC_n$ then  $Con(ZFC+IFS_{W\!Ab\,L(V_{\alpha})})$
\end{theorem}
\proof Let $\langle\theta_i\mid i<n\rangle$ be a list of formulas
defining weakly absolute sets of the form $V_{\alpha}.$   Let
$\big\{\varphi_{ij}(x)\mid i<n, j<m,\big\}$ be a set of formulas. It
is enough to show the consistency with $ZFC$ of 
$$\bigwedge\limits_{i<n,j<m}\exists
V_{\alpha_i}\Big[\theta_i(V_{\alpha_i})
\ \wedge\ \exists! P_{ij}(L(V_{\alpha})\models \varphi_{ij}(P_{ij}))\
\longrightarrow$$
$$\exists G(G\subseteq P_{ij}\ \wedge\ G\hbox{ is }L(V_{\alpha_i})\hbox{ generic}\Big]$$
Let $M$ be a countable transitive model of enough of $ZFC$ (including
$\mathop{\wedge}\limits_{i<n}ZF_{\theta_i}.\!$) Let
$\langle\alpha_0,\ldots,\alpha_{n-1}\rangle$ be the increasing sequence of ordinals such that
$$M\models \theta_i(V_{\alpha_i})$$ 
for $i<n.$ 
We define by induction on $(i,j)\in n\times m$ sets $G_{ij}.$ Suppose $P_{ij}$ is a partial order
definable in $L(V_{\alpha_i}^{M[\{G_{h,l}|h\leq
i,l<j\}]})$ by $\varphi_{ij}(x)$ and
there exists a $M[\{G_{h,l}|h\leq
i,l<j\}]$ generic subset of $P_{ij}$ not increasing 
$$V_{\alpha_i}^{M[\{G_{h,l}|h\leq
i,l<j\}]}$$ 
Then let $G_{ij}$ be such a $M[\{G_{h,l}|h\leq i,l<j\}]$
subset of $P_{ij}.$ (If not, let $G_{ij}=\emptyset.$) Let 
$$N=M[\{G_{ij}|i<n,j<m\}]$$
$N$ has the property that if $P_{ij}$ is a partial
order definable in $L(V_{\alpha_i})$ by $\varphi_{ij}(x)$ and $G$ is an $N$ generic subset
of $P_{ij}$ such that
$$\big|\big|V_{\alpha_i}^{N[G]}=V_{\alpha_i}^N\big|\big|^{(r.o.P)^N}=1$$
then an $L(V_{\alpha_i}^{N})$ generic subset of $P_{ij}$ exists in $N.$

\begin{theorem}
If $\langle M_n\mid n<\omega\rangle$ is a sequence of transitive models with
$M_n\models ZFC_n$ where $ZFC=\bigcup\limits_{n\in\omega}ZFC_n,$ then $Con(ZFC+IFS_{W\!Ab\,L[V_{\alpha}]}+IFS_{Ab\,L[r]})$
\end{theorem}
\proof Same as the last theorem except we start with a model of enough
of $ZFC+IFS_{Ab\,L[r]}.$

\begin{theorem}\label{Jech}
$V[G]$ has no functions $f:\kappa\rightarrow\kappa$ not in the ground
model if and only if $r.o.P$ is $(\kappa,\kappa)\!$-distributive.
\end{theorem}
\proof See [Jech1].

\begin{coro}
$IFS_{W\!Ab\,L(V_{\alpha})}$ is equivalent to the axiom schema of set
theory which says for every
weakly absolutely definable set of the form $V_{\alpha},$ for every
partial order 
$P$ definable in
$L(V_{\alpha}),$ if
$$(r.o.P)^V \hbox{ is }(\kappa,\kappa)\hbox{-distributive}$$
for each $\kappa$ such that for some $\beta<\alpha,\ \kappa\leq \big|V_{\beta}\big|,$  
then there exists an $L(V_{\alpha})$ generic subset $G$ of $P.$
\end{coro}

\begin{theorem}
$ZFC+IFS_{W\!Ab\ L(V_{\alpha})}\vdash\ CH$
\end{theorem}
\proof Let $P=$ the set of bijections from countable ordinals into
${\Bbb R}.$ Since $P$ is $\sigma$ closed, $\omega_1=\omega_1^{L({\Bbb
R})},$ and $P$ is a definable element of $L({\Bbb R}),$ there is an
$L({\Bbb R})$ generic subset of $P$ in $V.$ If $\alpha$ is an ordinal
less than $\omega_1$ and $r$ is a real, let $D_{\alpha}=\big\{p\in
P\mid \alpha\in dom\,p\big\}$ and $D_r=\big\{p\in
P\mid r\in ran\,p\big\}.$ For each $\alpha<\omega_1,$ $G\cap
D_{\alpha}\neq\emptyset$ and for each $r\in
{\Bbb R},$ $G\cap D_r\neq\emptyset,$ so $\bigcup G$ is a bijection from
$\omega_1$ to ${\Bbb R}.$

\vspace{.1in}

Perhaps the following is a better illustration of the kind of result
obtainable from $ZFC+IFS_{W\!Ab\ L(V_{\alpha})}.$

\begin{defi}
A Ramsey ultrafilter on $\omega$ is an Ultrafilter on $\omega$ such
that every coloring of $\omega$ with two colors has a homogenous set
in the ultrafilter.
\end{defi}

\begin{theorem}
$ZFC+IFS_{W\!Ab\ L(V_{\alpha})}\vdash $ there is a Ramsey ultrafilter
on $\omega.$
\end{theorem}
\proof Let $P$ be the partial order $(P(\omega),\sub^*)$ where
$P(\omega)$ is the power set of $\omega$ and $a\sub^* b$ means $a$ is
a subset of $b$ except for finitely many elements. $P$ is definable is
$L({\Bbb R})$ and is $\omega$ closed. The generic object is an Ramsey ultrafilter 
over $L({\Bbb R}),$ and since all colorings of $\omega$ are in $L({\Bbb R}),$ it is a
Ramsey ultrafilter over $V.$

\vspace{.1in}

\noindent One can argue that $IFS_{W\!Ab\ L(V_{\alpha})}$ is not a natural axiom
since among the definable sets $X$ with the property that $L(X)$ is
absolute when not increasing $X,$ why should you choose only
those of the form $V_{\alpha}?$ But it is natural in the sense it is a
way of forcing the universe as large as possible with respect to the
existence of generics by first fixing the
height of the models under consideration and then by fixing more and
more of their widths. In any case we should consider the
strengthenings of 
$IFS_{W\!Ab\,L(V_{\alpha})}$ defined below.

\begin{defi}
$x\in V$ is said to be weakly absolutely definable if for some
formula $\psi(x)$ which is provably $\Delta_1$ from $ZF,$ 
$$V\models \forall y(y\in x\ \leftrightarrow\ \psi(y))$$
\end{defi}

\begin{defi}
$IFS$ is the axiom schema of set theory which says for every
weakly absolutely definable set $X,$ for every partial order $P$ definable in
$L(X),$ if
$$\big|\big|X^{V[G]}=X^{V}\big|\big|^{(r.o.P)^V}=1$$
then there exists an $L(X)$ generic subset $G$ of $P.$
\end{defi}

\noindent If $X$ is an weakly absolutely definable set and $P$ is a partial ordering
definable in $L(X)$ such that 
$$\big|\big|X^{V[G]}=X^{V}\big|\big|^{(r.o.P)^V}=1$$
and if there is no $L(X)$ generic subset of $P$ in $V,$
we say that $V$ has a gap. $IFS$ says there are no gaps. 
The intuition that such gaps should not
occur in $V$ leads to the following: 

\begin{conj}
$ZFC+IFS$ is consistent.
\end{conj}

\noindent If $ZFC+IFS$ is consistent, then this means that it is
consistent that the universe is complete with respect to the
natural yardstick classes, (the classes of the form $L(X)$ where $X$
is weakly absolutely definable.) In my view, confirming the consistency of
$ZFC+IFS$ would be strong evidence that the universe of set theory
conforms to the axioms of $IFS.$ One reason for this opinion is that there is no apriori
reason for the consistency of $ZFC+IFS,$ so if $ZFC+IFS$ is
consistent, it seems that confirmimg its
consistency would involve some deep mathematics implying $IFS$
should be taken seriously.

\section{Formalizing the arguments in favor of $IFS_L$ and the other schemata}

\noindent In this section we try to formalize the vague notion that $IFS_L$ is a
natural closure condition on the universe, and that gaps in
general are esthetically undesirable. For simplicity we concentrate on $IFS_L.$

\begin{defi}
Let $T$ be a recursive theory in the language of set theory extending $ZFC.$  Let $P$
be a unary predicate.
If $\varphi$ is a formula of set theory
then $\varphi^*$ is $\varphi$ with all its quantifiers restricted to
$P,$ i.e., if $\exists x$ occurs in $\varphi$ then it is replaced by
$\exists x(P(x)\,\wedge\ldots)$ and $\forall x$ is replaced by
$\forall x(P(x)\ \rightarrow\ldots).$ 
The theory majorizing $T,$ $T',$ is the recursive theory in the
language $\big\{\varepsilon, P(x)\big\}$ such that
\begin{enumerate}
\item $\varphi\in T\ \rightarrow\ \varphi^*\in T'$
\item $P(x)\  is\  transitive\ \in T'$
\item $\forall x(x\in ORD\ \rightarrow\ P(x))\in T'$
\item $ZFC\subseteq T'$
\end{enumerate}

\noindent If $\theta(x)=\forall y(y\in x\leftrightarrow \psi(x))$ is a
formula defining a weakly absolutely 
definable set then the theory majorizing $T$ with respect to
$\theta(x)$ is $T'$ plus all the axioms of the form
$$\Big(\varphi_1\ \wedge\ \ldots\ \wedge\ \varphi_n\ \rightarrow\
\big(\psi(y)\leftrightarrow \exists z\psi_0(y,z)\leftrightarrow\forall
z\psi_1(y,z)\big)\Big)\ \longrightarrow\ \Big(\forall y(\psi(y)\ \rightarrow\ P(y))\Big)$$
where $\varphi_1,\ldots,\varphi_n\in ZF$ and $\psi_0(y,z)$ and
$\psi_1(y,z)$ are $\Delta_0$ formulas. 
\end{defi}

\begin{theorem}
Let $T$ be a recursive extension of $ZFC.$ Let
$T=\bigcup\limits_{n\in\omega}T_n$ where for some recursive function
$F,$ for each $n,$ $F(n)=T_n,$ a finite subset of $T$ and the $T_n$ are increasing. 
If there is a sequence $\langle M_n\mid n<\omega\rangle$ of
countable transitive models such that 
$$M_n\models T_n$$
then $T'+IFS_L$ ($T'$ is the theory majorizing $T\!$) is consistent
and there is a sequence $\langle N_n\mid n\in\omega\rangle$ of
countable transitive models such that $$N_n\models T_n'$$
where $T'=\bigcup\limits_{n\in\omega}T'_n$ and for some recursive
function $H,$ for each $n\in\omega,$ $H(n)=T'_n$ a
finite subset of $T'.$ 
\end{theorem}
\proof Let $IFS_L=\bigcup\limits_{n\in\omega}(IFS_L)_n$ where for each
$n\in\omega,\ (IFS_L)_n$ is finite. We can find a subsequence
$\langle N_n\mid n\in\omega\rangle$ of the $\langle
M_n\mid n\in\omega\rangle$ and $N_n$-generic sets $G_n$ such that 
$N_n[G_n]\models (IFS_L)_n,\ N_n\models T_n.$ Let $N_n[G_n]^*$
be the model in the language $\big\{\varepsilon,P(x)\big\}$ obtained
by letting the interpretation of $P(x)$ to be $N_n.$ Let $D$ be an
ultrafilter on $\omega.$ Then 
$$\prod N_n[G_n]^*/D$$
is a model for $T'+IFS_L.$

\begin{defi}
A theory extending $ZFC$ is $\omega-\!$complete if whenever $\varphi(x)$ is a
formula of set theory and if for each natural
number $n,$ 
$$T\vdash\varphi(n)$$
then $T\vdash \forall n\in\omega\varphi(n).$
\end{defi}

\begin{theorem}
Let $T$ be a recursive extension of $ZFC$ and suppose it has a
consistent, complete and $\omega-\!$complete extension $T^*.$ Then
$T'+IFS_L$ is 
consistent. 
\end{theorem}
\proof By reflection in $T^*,$ by its $\omega-\!$completeness and by the axiom of choice in $T^*,$
$$T^*\vdash\exists\langle N_n\mid n\in\omega\rangle$$
with the $\langle N_n\mid n\in\omega\rangle$ having  the same properties as in
the previous theorem. As in the previous theorem since $ZFC\subset T,$ $T^*\vdash Con(T'+IFS_L).$
Since $T^*$ is $\omega-\!$complete, (by the omitting
types theorem) it has an model $M$
with the standard set of integers. Since $M\models T^*,$ 
$$M\models Con(T'+IFS_L)$$
and as $Con(T'+IFS_L)$ is an arithmetical statement, it must really be
true. 

\vspace{.1in}

\noindent Certainly if the hypothesis of the theorem fails, then $T$
cannot be a suitable axiom system for set theory.

\begin{defi}
If $\theta(x)$ is a formula defining an weakly absolutely definable set,
then $IFS\restriction\theta(x)$ is $IFS$ restricted to the set defined by
$\theta(x),$ i.e., it says for all partial orders $P$ definable in
$L(X)$ were $X$ is defined by $\theta(x)$ such that
$$\big|\big|X^{V[G]}=X^{V}\big|\big|^{(r.o.P)^V}=1$$
there is an $L(X)$ generic
subset of $P.$
\end{defi}

\begin{theorem}
Let $\theta(x)$ be a formula defining an weakly absolutely definable set. 
Let $T$ be a recursive extension of $ZFC$ and suppose it has a
consistent, complete and $\omega-\!$complete extension $T^*.$ Then
$T_{\theta(x)}'+IFS\restriction\theta(x)$ is
consistent. 
\end{theorem}
\proof Same as above.

\begin{theorem}
Let $\theta(x)$ be a formula defining an weakly absolutely definable set. 
Let $T$ be a recursive extension of $ZFC+IFS\restriction\theta(x)$ and
suppose $T'_{\theta(x)}$ majorizes $T$ with respect to $\theta(x).$
Then $T_{\theta(x)}'\vdash IFS\restriction\theta(x).$  
\end{theorem}
\proof Working in $T'$ the generics in the inner model are still generic over $L(X)$
 since the inner model is a transitive class containing
all the ordinals.

\vspace{.1in}

\noindent The theorems in this section are meant as the formalization
of the notion that we can 'consistently imagine' a class model of
$ZFC$ not satisfying $IFS_L$ as being contained in a larger class
satisfying $ZFC+IFS_L,$ and that models of $ZFC$ not satisfying
$IFS$ have a gap.

\section{Conclusion}

\noindent These axiom schemata lead to many questions, among
them
\begin{enumerate}
\item Are there models of $IFS$ or $IFS_{W\!Ab\,L[V_{\alpha}]}$ which are
forcing extensions of $L$?
\item Are there similar natural schema's making the universe large,
but contradicting $IFS$ or
$IFS_{W\!Ab\,L[V_{\alpha}]}?$ 
\item What are the consequences for ordinary mathematics of these
axioms?
\end{enumerate}

\noindent The conventional view of the history of set theory says that Godel in
1938 proved that the consistency of $ZF$ implies the consistency of
$ZFC$ and of $ZFC+GCH,$ and that Cohen with the invention of forcing
proved that $Con(ZF)$ implies $Con(ZF+\neg AC)$ and $Con(ZFC+\neg
GCH)$ but from the point of view of $IFS_L$ a better way to state the history
would be to say that Godel
discovered $L$ and Cohen proved there are many generic sets over $L.$

I think confirming the consistency of $IFS$ with $ZFC$ would be a
vindication of the idea that generics over partial orders definable in
$L(X)$ with $X$ an weakly absolutely definable set exist, and
thereby put a crack in the armor surrounding the continuim hypothesis as
$ZFC+IFS\restriction{\Bbb R}\vdash CH.$ On the other hand, if
$ZFC+IFS$ is not consistent, it would show the universe must have some
gaps, i.e., incomplete with respect to some concrete set, an
esthetically unpleasing result. 
It is ironic that although mathematics and especially mathematical
logic is an art noted for its precise and formalized reasoning, 
it seems that in order to solve problems at the frontiers of logic's
foundations we must tackle questions of an esthetic nature.

\pagebreak

\begin{center}
REFERENCES
\end{center}

\noindent 1. [CK] C. C. Chang and J. Keisler, {\em Model Theory},
North Holland Publishing Co.

\vspace{.1in}

\noindent 2. [Jech1] T. Jech, {\em Multiple Forcing}, Cambridge
University Press.

\vspace{.1in}

\noindent 3. [Jech2] T. Jech, {\em Set Theory}, Academic Press.

\end{document}